\documentclass{article}
\usepackage{amsmath}
\usepackage{amsthm}
\usepackage{graphicx}
\usepackage{algorithm}
\usepackage{algorithmic}
\usepackage{amssymb}
\usepackage{geometry}
\usepackage{mathtools}
\usepackage{color}

\newtheorem{theorem}{Theorem}
\newtheorem{definition}{Definition}
\newtheorem{remark}{Remark}
\newtheorem{proposition}{Proposition}
\newtheorem{example}{Example}
\newtheorem{lemma}{Lemma}

\title{A Systematic Review on Hermite-Hadamard Inequality: Theory and Applications}
\author{Ohud Almutairi$^{1}$ and Adem K{\i}l{\i}\c{c}man$^{2}$	}

\date{\today}

\begin{document}
	\maketitle
	
	\begin{center} $^{1}$Department of Mathematics, University of Hafr Al Batin, Hafr Al-Batin 31991,Sudia Arabia. OhudbAlmutairi@gmail.com \\ $^{2}$Department of Mathematics and Statistics, University Putra Malaysia, 43400 UPM, Serdang, Selangor, Malaysia, akilic@upm.edu.my \end{center}
	\noindent \textbf{Abstract}.
Inequalities play important roles not only in mathematics, but also in other fields, such as economics and engineering. Even though many results are published on Hermite-Hadamard (H-H) type inequalities, new researcher to this fields often found it difficult to understand them. Thus, some important discoverers, such as the formulations of H-H type inequalities via various classes of convexity, through differentiable mappings and for fractional integrals, are presented. Some well-known examples from previous literature are used as illustrations.

\section{Introduction}
Mathematical inequalities play a key role in understanding a range of problems in various fields of mathematics. Among the most celebrated ones is Hermite-Hadamard (H-H) inequality, which made a great impact not only in mathematics, but also in other related disciplines.
As mentioned by Mitrinovi\'{c} and Lackovi\'{c} \cite{mitrinovic1985hermite}, this inequality was first appeared in the literature through the effort of Hadamard \cite{hadamard1893etude}; however, the result was first discovered by Hermite \cite{hermite1883deux}. Following this fact, many researchers referred the result as the H-H inequality. This inequality was stated in the monograph of \cite{dragomir2004selected} as the first fundamental result for convex functions defined in the interval of real numbers with a natural geometrical interpretation that can be applied to investigate a variety of problems. Inequalities play important roles in understanding many mathematical concepts, such as probability theory, numerical integration and integral operator theory. Throughout the last century, the H-H type inequalities have been considered among the fast growing fields in mathematical analysis, through which vast problems in engineering, economics and physics have been studied \cite{bullen2003handbook,dragomir2004selected,wang2018fractional}. Due to the enormous importance of these inequalities, many extensions, refinements and generalizations of their related types have been equally investigated \cite{Ohudalmutairil,almutairi2019new,bin2019new,duc2020convexity}.

Therefore, the H-H type inequalities, by which many results are studied, play important roles in the theory of convex functions. 
The convexities, along with many types of their generalizations, can be applied in different fields of sciences \cite{liu2020convex,almutairi2020new}, through which many generalizations of H-H inequality for a variant types of convexities have been studied. Other extensions of H-H inequality include the formulation of problems related to fractional calculus, a branch of calculus dealing with derivatives and integrals of non-integer order \cite {dragomir2019hermite,dahmani2020integral,almutairi2021generalized}.

This paper is aimed at introducing the H-H inequality to a new researcher in the field. Thus, we present basic facts on some integral inequalities, fractional inequalities of H-H type and their constructions via various convexity classes. Some important theorems associated with these inequalities are also discussed, along with some well-known examples, to ease the beginners’ understanding of the basic concepts of these inequalities. Even though the information presented in this review article can be found in separate studies of inequalities, obtaining a single work combining these results remains elusive. 

\section{preliminaries}\label{cp1}
In the following, we will give some necessary definitions
and mathematical preliminaries of fractional calculus
theory which are used further in this paper. For more
details, one can consult \cite{udriste2013convex} and \cite{ullah2019note}. The concept of a convex function was first introduced to elementary calculus when discussing necessary conditions for a minimum or maximum value of a differentiable function. The convex function was later recognized as an active area of study by  \cite{jensen1905om}.
In modern studies, a convex function is considered as a link between analysis and geometry, which makes it a powerful tool for solving many practical problems.

\begin{definition}\label{definitionconvexity} \cite{niculescu2006convex} Let $V$ be an interval in $\mathbb{R}$. A function \(\mathcal{G}: V \rightarrow \mathbb{R}\) is said to be convex if
	\begin{equation}
	\mathcal{G}(\zeta m_1+(1-\zeta) m_2) \leq\zeta \mathcal{G}(m_1)+(1-\zeta) \mathcal{G}(m_2) \label{definitionconvex}
	\end{equation}
	holds for all \(m_1, m_2 \in V\) and \(\zeta \in[0,1]\).	
\end{definition}

If inequality (\ref{definitionconvex}) strictly holds for any distinct points $m_1$ and $m_2$, where  \(\zeta \in(0,1) \), then the function is said to be a strictly convex. Meanwhile, a function \(-\mathcal{G}\) is convex (strictly convex), then \(\mathcal{G}\) is concave (strictly concave).

Geometrically, a function $\mathcal{G}$ is convex given that the line segment joining any two points on the graph lies above (or on) the graph. Meanwhile, if the line segment connecting the two points is below (or on) the graph, the function is concave.

\begin{example} Given a function \(\mathcal{G}: V\subseteq \mathbb{R}\rightarrow \mathbb{R}\) for any $m \in \mathbb{R}$, we have the following examples.
	\begin{itemize}
		\item [i.] $\mathcal{G}(m)=c_1 m+c_2$, where \(c_1, c_2 \in \mathbb{R}.\) The function $\mathcal{G}(m)$ is both concave and convex on \((-\infty, \infty)\). Thus, it is  referred to as an affine.
		\item[ii. ] The functions $\mathcal{G}(m)=m^{2}$ and $\mathcal{G}(m)=e^{m}$ are both convex functions on \(\mathbb{R}\).
		\item[iii.] $\mathcal{G}(m)=\ln m$ is a concave function on $\mathbb{R}_{+}=[0,\infty)$.
	\end{itemize}
\end{example}


The theory of convexity deals with large classes, such as generalized convex function on fractal sets, Godunova-Levin, $s$-convex and preinvex functions. These, termed as the generalization of convexity, play important roles in optimization theory and mathematical programming.Therefore, we give basic definitions of different classes of convex functions.

The definition of generalized convex functions on fractal sets \(\mathbb{R}^{\alpha}(0<\alpha\leq1)\) is given by  \cite{mo2014generalized} as follows.

\begin{definition}
	Let \(\mathcal{G}: V \subseteq \mathbb{R} \rightarrow \mathbb{R}^{\alpha} .\) For any \(m_{1}, m_{2} \in V\) and \(\zeta \in[0,1],\) if the following inequality
	\begin{equation*}
	\mathcal{G}\left(\zeta m_1+(1-\zeta) m_2\right) \leq \zeta^{\alpha} \mathcal{G}\left(m_1\right)+(1-\zeta)^{\alpha} \mathcal{G}\left(m_2\right)
	\end{equation*}
	holds, then \(\mathcal{G}\) is called a generalized convex function on $V$.
\end{definition}

The space Godunova-Levin function, denoted by $Q(V)$, was introduced by \cite{godunova1985neravenstva}. They noted that both the positive monotone and positive convex functions are belonged to $Q(V)$. Due to the importance of this function, we present it as follows.

\begin{definition}\cite{mitrinovic2013classical}
	A non-negative function $\mathcal{G}: V \subseteq \mathbb{R} \rightarrow \mathbb{R}$ is called Godunova-Levin function (denoted by
	$\mathcal{G} \in Q(V)$) if
	\begin{equation}
	\mathcal{G}(\zeta m_1+(1-\zeta) m_2) \leq \frac{1}{\zeta} \mathcal{G}(m_1)+\frac{1}{1-\zeta} \mathcal{G}(m_2)\label{GodunovaLevin}
	\end{equation}
	holds, for all \(m_1, m_2 \in V\) and \(\zeta \in(0,1)\).
\end{definition}

\begin{example}\cite{dragomir1995some}
	For \(x \in[m_1, m_2],\) the function
	
	\begin{equation*}
	\mathcal{G}(x)=\left\{\begin{array}{ll}{1,} & {m_1 \leq x<\frac{m_1+m_2}{2}} \\ {4,} & {x=\frac{m_1+m_2}{2}} \\ {1,} & {\frac{m_1+m_2}{2}<x \leq m_2}\end{array}\right.
	\end{equation*}	
	is in the class \(Q(V)\).
\end{example}


Godunova-Levin function was restricted to a space called $P(V)$ contained in $Q(V)$. This class is defined by \cite{dragomir1995some} as follows.

\begin{definition}
	A non-negative function $\mathcal{G}: V \subseteq \mathbb{R}\rightarrow \mathbb{R}$ is called P-function (denoted by $\mathcal{G} \in P(V)$) if
	\begin{equation*}
	\mathcal{G}(\zeta m_1+(1-\zeta) m_2) \leq \mathcal{G}(m_1)+\mathcal{G}(m_2)
	\end{equation*}
	holds, for all \(m_1, m_2 \in V\) and \(\zeta \in[0,1]\).	
\end{definition}
Therefore, all non-negative monotone and convex functions are contained in \(P(V)\).

For other results of Godunova-Levin and $P$-functions, see \cite{radulescu2009godunova}, \cite{fang2014p} and \cite{kadakal2017some}.

The definition of $s$-convex function in the second sense 
or $s$-Breckner convex is given as follows.
\begin{definition}\label{ch1-defi-sconvex}\cite{breckner1978stetigkeitsaussagen}
	A function \(\mathcal{G}: [0,\infty) \rightarrow \mathbb{R}\) is said to be \(s\)-convex in second sense (denoted by $\mathcal{G} \in K_{s}^{2}$), if
	\begin{equation}
	\mathcal{G}(\zeta_1 m_1+\zeta_2 m_2) \leq \zeta_1^{s} \mathcal{G}(m_1)+\zeta_2^{s} \mathcal{G}(m_2)\label{ssconvex}
	\end{equation}
	holds, for all \(m_1, m_2 \in [0,\infty)\), $\zeta_1,\zeta_2\geq0$, $\zeta_1+\zeta_2=1$ and $0<s \leq 1$.	
\end{definition}

Choosing $s=1$ reduces $s$-convexity in second sense to the classical convex function on \([0, \infty)\).

The following property that is connected to \(s\)-convex function in the second sense is given bellow.

\begin{theorem}\cite{hudzik1994some}
	If \(\mathcal{G} \in K_{s}^{2},\) then \(\mathcal{G}\) is non-negative on \([0, \infty) .\)
\end{theorem}
For some properties of $s$-convexity in second sense, see the references \cite{dragomir1999hadamard,du2017generalization,usta2018generalization,gozpinar2019some}.

\cite{hudzik1994some} present the example of \(s\)-convex function in the second sense as follows.

\begin{example}
	Let 0 \(<s<1\) and \(c_1, c_2, c_3 \in \mathbb{R} .\) Defining

	\begin{equation*}
	\mathcal{G}(m)=\left\{\begin{array}{ll}{c_1,} & {m=0} \\ {c_2 m^{s}+c_3,} & {m>0}\end{array}\right.
	\end{equation*}	
	for  \(m\in \mathbb{R}_{+}\), we have
	\begin{itemize}
		\item [i.] If \(c_2 \geq 0\) and \(0 \leq c_3 \leq c_1,\) then \(\mathcal{G} \in K_{s}^{2}\),
		\item 	[ii.]If \(c_2>0\) and \(c_3<0,\) then \(\mathcal{G} \notin K_{s}^{2}\).
	\end{itemize}
\end{example}

As Hudzik and Maligranda mentioned that the condition
\(\zeta_1+\zeta_2=1\) in definition \ref{ch1-defi-sconvex} can be replaced by \(\zeta_1+\zeta_2 \leq1\), equivalently.

\begin{theorem}\cite{hudzik1994some}
	Suppose that \(\mathcal{G} \in K_{s}^{2}.\) The inequality
	(\ref{ssconvex}) holds for all \(c_1, c_2 \in \mathbb{R}_{+}\) and \(\zeta_1, \zeta_2 \geq 0\) with \(\zeta_1+\zeta_2 \leq 1\) iff \(\mathcal{G}(0)=0\).
\end{theorem}

The geometric description of \(s\)-convex curve, given in the definition below, was clearly explained in \cite{pinheiro2007exploring}.

\begin{definition}
	A function $\mathcal{G}:V\subseteq \mathbb{R}\to \mathbb{R}$ is called an $s$-convex in the second sense for \(0<s<1\), if the graph of the function is below a bent chord $L$ that is between any two points. This means that, for every compact
	interval \(W \subset V,\) the following inequality	
	\begin{equation*}
	\sup _{W}(L-\mathcal{G}) \geq \sup _{\partial W}(L-\mathcal{G})
	\end{equation*}
	holds, with boundary \(\partial W.\) 
\end{definition}

The $s$-convex function of second sense can be referred
as the limiting curve. This differentiates the curves of $s$-convex
in second sense from others which are not. Following this, Pinheiro determines the affects of the choice
of $s$ on the limiting curve. For further results on $s$-convex function in the second sense, we refer the reader to  \cite{dragomir2000jensen}, \cite{alomari2008hadamard} and \cite{dragomir2016integral}.

The definition of the generalized $s$-convex function  on fractal sets is given as follows.
\begin{definition}\cite{mo2014generalized} A function $\mathcal{G} :V \subseteq \mathbb{R}_{+} \rightarrow \mathbb{R}^{\alpha}$ is a generalized $s$-convex in the second sense on fractal sets if
	\begin{equation}
	\mathcal{G}\left(\zeta_1 m_1+\zeta_2 m_2\right) \leq \zeta_1^{\alpha s} \mathcal{G}(m_1)+\zeta_2^{\alpha s} \mathcal{G} (m_2)\label{sconvexfractalset}
	\end{equation}
	holds, for all \(m_1, m_2 \in V\), $0<s\leq 1$, $\zeta_1,\zeta_2\geq0$ and $\zeta_1+\zeta_2=1$. This class of function is denoted by $GK_s^{2}$.
\end{definition}

The generalized $s$-convex function in the second sense becomes $s$-convex function when $\alpha=1$. 

One should note that the following theorems along with example can be found in \cite{mo2014generalized}.

\begin{theorem}
	Let \(\mathcal{G} \in G K_{s}^{2} .\) Inequality (\ref{sconvexfractalset}) holds for all \(m_1, m_2 \in \mathbb{R}_{+}\) and \(\zeta_{1}, \zeta_{2} \geq 0\)	with \(\zeta_{1}+\zeta_{2}<1\) iff \(\mathcal{G}(0)=0^{\alpha}\).
\end{theorem}

\begin{theorem}
	Let \(0<s<1\). If \(\mathcal{G} \in G K_{s}^{2},\) then \(\mathcal{G}\) is non-negative on \([0,+\infty)\).
\end{theorem}

\begin{theorem}
	Let \(0<s_{1} \leq s_{2} \leq 1 .\) If \(\mathcal{G} \in G K_{s_{2}}^{2}\) and $\mathcal{G}(0)=0^{\alpha}$, then
	\(\mathcal{G} \in G K_{s_{1}}^{2}\).
\end{theorem}

Considering the properties of the generalized $s$-convex in the second sense, we present the following example.

\begin{example}
	Let \(0<s<1,\) and \(a_1^{\alpha}, a_2^{\alpha}, a_3^{\alpha} \in \mathbb{R}^{\alpha} .\) For \(m \in \mathbb{R}_{+},\) we define
	
	\begin{equation*}
	\mathcal{G}(m)=\left\{\begin{array}{ll}{a_1^{\alpha},} & {m=0} \\ {a_2^{\alpha} m^{s \alpha}+a_3^{\alpha},} & {m>0}.\end{array}\right.
	\end{equation*}
	Thus, we have the following:
	\begin{itemize}
		\item [i.] If \(a_2^{\alpha} \geq 0^{\alpha}\) and \(0^{\alpha} \leq a_3^{\alpha} \leq a_1^{\alpha},\) then \(\mathcal{G} \in G K_{s}^{2}\),
		\item [ii.] If \(a_2^{\alpha}>0^{\alpha}\) and \(a_3^{\alpha}<0^{\alpha},\) then \(\mathcal{G} \notin G K_{s}^{2}\).
	\end{itemize}	
\end{example}

. 

For more results related to the generalized $s$-convex in the second sense on fractal sets , the interested reader is directed to \cite{kiliccman2015notions} and \cite{budak2016generalized}.

In order to unify the concepts of Godunova-Levin and P-functions, \cite{dragomir2016integral}
introduced $s$-Godunova-Levin as follows.

\begin{definition} A function $\mathcal{G}: V\subseteq \mathbb{R} \rightarrow[0, \infty)$ is said to be $s$-Godunova-Levin, (denoted by $\mathcal{G} \in Q_{s_{2}}(V)$), if
	\begin{equation}
	\mathcal{G}(\zeta m_1+(1-\zeta) m_2) \leq \frac{1}{\zeta^{s}} \mathcal{G}(m_1)+\frac{1}{(1-\zeta)^{s}} \mathcal{G}(m_2)\label{sGodunovaLevine}
	\end{equation}	
	holds, for all \(m_1, m_2 \in V\), $\zeta\in (0,1)$ and $0\leq s \leq 1$.	
\end{definition}

Choosing $s=1$ reduces $s$-Godunova-Levin function to the the class of Godunova-Levin.
Also, taking $s=0$ we have the class of \(P\)-function. Thus, we have the following,
$P(V)=Q_{0}(V) \subseteq Q_{s_{2}}(V) \subseteq Q_{1}(V)=Q(V)$. For more results on $s$-Godunova-Levin of convexity, we refer the reader to \cite{noor2014fractional} and \cite{kashuri2018hermite}.

Preinvex functions are among the most important classes of generalized convex functions. This concept, playing important roles in many disciplines, was proposed by \cite{ben1986invexity}. Since then, a preinvex function has become an active area of study. 

\begin{definition}\cite{hanson1981sufficiency}
	A set \(V \subseteq \mathbb{R}\) is called an invex if there exists a function $\eta: V \times V \rightarrow \mathbb{R}$ such that
	\begin{equation*}
	m_1+\zeta \eta(m_2, m_1) \in V
	\end{equation*}
	holds, for all $m_1, m_2 \in V$ and $\zeta \in[0,1]$.	
\end{definition}
The invex set $V$ can also referred to as an $\eta$-connected set.

\begin{definition}\cite{ben1986invexity}\label{defpreinvexity}
	Suppose that $V\subseteq\mathbb{R}$ is an invex set with respect to \(\eta: V \times V \rightarrow \mathbb{R}\). A function \(\mathcal{G}: V \rightarrow \mathbb{R}\) is called preinvex with respect to \(\eta,\) if
	\begin{equation}
	\mathcal{G}(m_1+\zeta \eta(m_2, m_1)) \leq(1-\zeta) \mathcal{G}(m_1)+\zeta \mathcal{G}(m_2)\label{preinvexdefinition}
	\end{equation}
	holds, for all $m_1, m_2\in V$ and $\zeta\in [0,1]$.
\end{definition}

Further generalizations can be found in \cite{weir1988pre,jue2010hadamard,meftah2017hadamard,meftah2019fractional}.

Fractional calculus, whose applications can be found in many disciplines including economics, life and physical sciences, as well as engineering, can be considered as one of the modern branches of mathematics \cite{cafagna2007fractional,machado2011recent,yang2019general,hilfer2019mathematical}. Many problems of interests from these fields can be analyzed through fractional integrals, which can also be regarded as an interesting sub-discipline of fractional calculus. Some of the applications of integral calculus can be seen in the following papers [5,6,7,8,9,10], through which problems in physics, chemistry, and population dynamics were studied. The fractional integrals were extended to include the H–H inequality \cite{baleanu2011fractional}, \cite{nigmatullin1992fractional}, \cite{kilbas2006theory}, \cite{diethelm2010analysis}, \cite{tarasov2011fractional}, \cite{malinowska2015advanced}, \cite{salati2019direct}. Now, we recall some basic definitions of fractional integrals as follows

\begin{definition}
	Let $\mathcal{G} \in L_{1}[m_1, m_2]$. The left and right sides Riemann-Liouville integrals denoted by $J_{m_1^+}^{\lambda} \mathcal{G}$ and $J_{m_2^-}^{\lambda} \mathcal{G} $ of order $\lambda\in \mathbb{R}_{+}$  are defined by
	\begin{flalign*}
	\quad J_{m_1^+}^{\lambda} \mathcal{G}(x)=\frac{1}{\Gamma(\lambda)} \int_{m_1}^{x}(x-\gamma)^{\lambda-1} \mathcal{G}(\gamma) d \gamma, \quad x>m_1
	\end{flalign*}
	and
	\begin{flalign*}
	\qquad J_{m_2^-}^{\lambda} \mathcal{G}(x)=\frac{1}{\Gamma(\lambda)} \int_{x}^{m_2}(\gamma-x)^{\lambda-1} \mathcal{G}(\gamma) d \gamma, \quad x<m_2,
	\end{flalign*}
	respectively. 
\end{definition}

If $\lambda=1$ in the above equalities, we get the classical integral.

One should note that the Hadamard fractional integrals differ from those of the Riemann-Liouville since in the former the logarithmic functions of arbitrary exponents are included in the kernels of the integrals. Therefore, the Hadamard fractional integrals are defined as follows.

\begin{definition}\cite{samko1993fractional}
	Let \(\lambda>0\) with \(m-1<\lambda \leq m, m \in \mathbb{N},\) and \(m_1<x<m_2.\) The left and right sides Hadamard fractional integrals denoted by $H_{m_1^+}^{\lambda}\mathcal{G}(x)$ and $H_{m_2^-}^{\lambda} \mathcal{G}(x)$ of order \(\lambda\) of a function \(\mathcal{G}\) are given as
	\begin{equation*}
	H_{m_1^+}^{\lambda} \mathcal{G}(x)=\frac{1}{\Gamma(\lambda)} \int_{\mathrm{m_1}}^{x}\left(\ln \frac{x}{\gamma}\right)^{\lambda-1} \frac{\mathcal{G}(\gamma)}{\gamma} d\gamma
	\end{equation*}
	and
	\begin{equation*}
	H_{m_2^-}^{\lambda} \mathcal{G}(x)=\frac{1}{\Gamma(\lambda)} \int_{x}^{\mathrm{m_2}}\left(\ln \frac{\gamma}{x}\right)^{\lambda-1} \frac{\mathcal{G}(\gamma)}{\gamma} d \gamma,
	\end{equation*}
	respectively.
\end{definition}

\cite{anatoly2001hadamard}, \cite{butzer2002compositions,butzer2002fractional} and \cite{kilbas2006theory} provide useful background and properties of Hadamard fractional integrals.

The following proposition is related to the Hadamard integrals.

\begin{proposition}\cite{kilbas2006theory}
	If \(\lambda>0\) and $0<m_1<m_2<\infty$, the following relations hold:
	\begin{equation*}
	\bigg(H_{m_1^+}^{\lambda}\bigg(\log \frac{\gamma}{m_1}\bigg)^{\beta-1}\bigg)(x)=\frac{\Gamma(\beta)}{\Gamma(\beta+\lambda)}\bigg(\log \frac{x}{m_1}\bigg)^{\beta+\lambda-1}
	\end{equation*}
	and
	\begin{equation*}
	\bigg(H_{m_2^-}^{\lambda}\bigg(\log \frac{m_2}{\gamma}\bigg)^{\beta-1}\bigg)(x)=\frac{\Gamma(\beta)}{\Gamma(\beta+\lambda)}\bigg(\log \frac{m_2}{x}\bigg)^{\beta+\lambda-1}.
	\end{equation*}
\end{proposition}

The Riemann-Liouville fractional integrals, along with the Hadamard’s fractional integrals, are generalized through the recent work of \cite{katugampola2015mellin}. These two integrals were combined and given in a single form. The following definition \cite{katugampola2015mellin} modifies the old version \cite{katugampola2011new} for Katugampola fractional integrals. 

\begin{definition}
	Let \([m_1, m_2] \subset \mathbb{R}\) be a finite interval. The left and right-sided Katugampola
	fractional integrals of order $\lambda>0$ for \(\mathcal{G} \in X_{c}^{p}(m_1, m_2)\) are defined by
	\begin{equation*}
	^\rho I_{m_1^+}^{\lambda} \mathcal{G}(x)=\frac{\rho^{1-\lambda}}{\Gamma(\lambda)} \int_{\mathrm{m_1}}^{x} \frac{\gamma^{\rho-1}}{\left(x^{\rho}-\gamma^{\rho}\right)^{1-\lambda}}\mathcal{G}(\gamma) d \gamma 
	\end{equation*}
	and
	\begin{equation*}
	\quad^{\rho} I_{m_2^-}^{\lambda} \mathcal{G}(x)=\frac{\rho^{1-\lambda}}{\Gamma(\lambda)} \int_{x}^{m_2} \frac{\gamma^{\rho-1}}{\left(\gamma^{\rho}-x^{\rho}\right)^{1-\lambda}} \mathcal{G}(\gamma) d\gamma,
	\end{equation*}
	with $m_1<x<m_2$ and $\rho>0$. 
\end{definition}

Following this, the space \(X_{c}^{p}(m_1, m_2)(c \in \mathbb{R}, 1 \leq p \leq \infty)\) is introduced as follows.
\begin{definition}\cite{anatoly2001hadamard} Let the space \(X_{c}^{p}(m_1, m_2)(c \in \mathbb{R}, 1 \leq p \leq \infty)\) of those complex-valued Lebesgue measurable functions \(\mathcal{G}\) on \([m_1, m_2]\) for which \(\|\mathcal{G}\|_{X_{c}^{p}}<\infty,\) where the norm is defined by
	\begin{equation*}
	\|\mathcal{G}\|_{X_{c}^{p}}=\left(\int_{m_1}^{m_2}\left|\zeta^{c} \mathcal{G}(\zeta)\right|^{p} \frac{d \zeta}{\zeta}\right)^{1 / p}<\infty \quad(1 \leq p < \infty, c \in \mathbb{R})
	\end{equation*}
	and for the case \(p=\infty\)
	\begin{equation*}
	\|\mathcal{G}\|_{X_{c}^{\infty}}=\operatorname{ess} \sup _{m_1 \leq \zeta \leq m_2}\left(\zeta^{c}|\mathcal{G}(\zeta)|\right) \quad(c \in \mathbb{R}),
	\end{equation*}
	where ess sup \(|\mathcal{G}(\zeta)|\) stands for the essential maximum of \(|\mathcal{G}(\zeta)|\). 
\end{definition}

If \(c=1 / p,\) \(X_{c}^{p}(m_1, m_2)\) reduces to \(L_{p}(m_1, m_2),\) the $p$-integrable function.

Important references on Katugampola fractional integrals and their applications are suggested for further reading \cite{butkovskii2013fractional, gaboury2013some, richard2014fractional, katugampola2014new}.

The relations among Katugampola fractional integrals, Riemann-Liouville integrals and Hadamard integrals are given in the next theorem. The left-sided version of the relation is considered here for its simplicity since similar results also exist for the right-sided operators.

\begin{theorem}\cite{katugampola2014new}
	Let \(\lambda>0\) and \(\rho>0 .\) Then for \(x>m_1\), we have
	\begin{itemize}
		\item [i.] 	$\lim_{\rho \rightarrow 1}$ $^{\rho}I_{m_1^+}^{\lambda} \mathcal{G}(x)=J_{m_1^+}^{\lambda} \mathcal{G}(x)$,
		\item [ii.] $\lim _{\rho \rightarrow 0^{+}}$   $^{\rho}I_{m_1^+}^{\lambda} \mathcal{G}(x)=H_{m_1^+}^{\lambda} \mathcal{G}(x)$.
	\end{itemize}	
\end{theorem}

\begin{remark}
	One should note that, while (i) is concerned with the Riemann-Liouville operators, (ii) is related to the Hadamard operators. 
\end{remark}


The definitions of the conformable fractional derivative and integral were given in \cite{khalil2014new}, and we present them as follows.

\begin{definition}
Let \(\mathcal{G}:[0, \infty) \rightarrow \mathbb{R}\), the conformable fractional derivative of \(\mathcal{G}\) of order \(\alpha\) is defined as

\begin{equation}
D_{\alpha}(\mathcal{G})(r)=\lim _{b \rightarrow 0} \frac{\mathcal{G}\left(r+b r^{1-\alpha}\right)-\mathcal{G}(r)}{b}
\end{equation}
where \(\mathcal{G}\) is said to be \(\alpha\)-differentiable at \(r\) if \(D_{\alpha}(\mathcal{G})(r)\) exists. In particular, \(D_{\alpha}(\mathcal{G})(0)\) is
defined as follows
\begin{equation*}
D_{\alpha}(\mathcal{G})(0)=\lim _{r \rightarrow 0^{+}} D_{\alpha}(\mathcal{G})(r)
\end{equation*}
and we use \(\mathcal{G}^{\alpha}(r)\) or \(\left(\mathrm{d}_{\alpha} / \mathrm{d}_{\alpha} r\right)(\mathcal{G})\) to denote \(D_{\alpha}(\mathcal{G})(r)\)
\end{definition}

\begin{definition}
Let \(\alpha \in(0,1]\) and \(0 \leq m_1<m_2\). A function
\(\mathcal{G}:[m_1, m_2] \rightarrow \mathbb{R}\) is \(\alpha\)-fractional integrable on \([m_1, m_2]\)
if the integral
\begin{equation*}
\int_{m_1}^{m_2} \mathcal{G}(x) d_{\alpha} x:=\int_{m_1}^{m_2} \mathcal{G}(x) x^{\alpha-1} d x
\end{equation*}
exists and is finite. All \(\alpha\)-fractional integrable on
\([m_1, m_2]\) is indicated by \(L_{\alpha}^{1}([m_1, m_2])\)
\end{definition}

\begin{remark}
	\begin{equation*}
	I_{\alpha}^{m_1}(\mathcal{G})(\gamma)=I_{1}^{m_1}\left(\gamma^{\alpha-1} \mathcal{G}\right)=\int_{m_1}^{t} \frac{\mathcal{G}(x)}{x^{1-\alpha}} d x
	\end{equation*}
	where the integral is the usual Riemann improper
	integral, and \(\alpha \in(0,1] .\)
\end{remark}

\begin{theorem}
 Let \(\alpha \in(0,1]\) and \(\mathcal{G}:[m_1, m_2] \rightarrow \mathbb{R}\) be
	a continuous on \([m_1, m_2]\) with \(0 \leq m_1<m_2\). Then,
	\begin{equation*}
	\left|I_{\alpha}^{m_1}(\mathcal{G})(x)\right| \leq I_{\alpha}^{m_1}|\mathcal{G}|(x)
	\end{equation*}
\end{theorem}
For more results on conformable integral operators, we refer the interested reader to \cite{khurshid2019conformable} and \cite{iqbal2020some}.


The H\"older integral inequality plays an important role in both pure and applied sciences. Other areas of applying this inequality include the theory of convexity,  which can be considered as one of the active and fast growing fields of studies in mathematical sciences. Thus, the H\"older's integral inequality is described in the following theorem.

\begin{theorem}\cite{mitrinovic1970analytic}
	Suppose that $p>1$ and \(1 / p+1 / q=1\). If \(\mathcal{G}\) and $\mathcal{K}$ are real functions on \([m_1, m_2]\) such that \(|\mathcal{G}|^{p}\) and \(|\mathcal{K}|^{q}\) are integrable functions on \([m_1, m_2],\) then
	\begin{equation*}
	\int_{m_1}^{m_2}|\mathcal{G}(x) \mathcal{K}(x)| d x \leq\left(\int_{m_1}^{m_2}|\mathcal{G}(x)|^{p} d x\right)^{\frac{1}{p}}\left(\int_{m_1}^{m_2}|\mathcal{K}(x)|^{q} d x\right)^{\frac{1}{q}}
	\end{equation*}
	holds.	
\end{theorem}


The other version of H\"older integral inequality is called the power-mean integral, which is given in the following theorem.
\begin{theorem}\cite{mitrinovic2013classical} 
	Suppose that \(q \geq 1 .\) Let $\mathcal{G}$ and $\mathcal{K}$ be real mappings on $[m_1, m_2]$. If $|\mathcal{G}|$ and $|\mathcal{G}||\mathcal{K}|^{q}$ are integrable functions in the given interval,  then
	\begin{equation*}
	\int_{m_1}^{m_2}|\mathcal{G}(x) \mathcal{K}(x)| d x \leq\left(\int_{m_1}^{m_2}|\mathcal{G}(x)| d x\right)^{1-\frac{1}{q}}\left(\int_{m_1}^{m_2}|\mathcal{G}(x)||\mathcal{K}(x)|^{q} d x\right)^{\frac{1}{q}}
	\end{equation*}
	holds.
\end{theorem}

\section{Hermite-Hadamard inequality}
H-H inequality plays a vital role in the theory of convexity.
This inequality estimates the integral average of any convex functions through the midpoint and trapezoidal formula of a given domain. While the midpoint formula estimates the integral from the left, the trapezoidal formula estimates it from the right.
More precisely, the classical H-H inequality is considered as follows.

\begin{theorem}\cite{dragomir2004selected}
	Let $\mathcal{G}:[m_1,m_2]\subseteq\mathbb{R}\to \mathbb{R}$ be a convex function on $[m_1,m_2]$ with $m_1<m_2$, then
	
	\begin{equation}
	(m_2-m_1)\mathcal{G}\left(\frac{m_1+m_2}{2}\right) \leq \int_{m_1}^{m_2} \mathcal{G}(x) d x \leq (m_2-m_1)\frac{\mathcal{G}(m_1)+\mathcal{G}(m_2)}{2}\label{hh}
	\end{equation}
	holds.
\end{theorem}

The proof of inequality (\ref{hh}) is provided here for simplicity. Though the proof of the theorem exists, this is the first time (\ref{hh}) is proved using a similar technique reported in \cite{sarikaya2013hermite}.

\begin{proof}
	Let $\mathcal{G}$ be a convex function on the interval $[m_1,m_2]$. Taking $\zeta=\frac{1}{2}$ in inequality (\ref {definitionconvex}) for $x,y\in [m_1,m_2]$, we have
	\begin{equation}
	\mathcal{G}\left(\frac{x+y}{2}\right)\leq \frac{\mathcal{G}(x)+\mathcal{G}(y)}{2}.\label{ch1proofhh1}
	\end{equation}
	Substituting $x=\zeta m_1+(1-\zeta)m_2$ and $y=(1-\zeta)m_1+\zeta m_2$ in (\ref{ch1proofhh1}), we get
	\begin{equation}
	2\mathcal{G}\left(\frac{m_1+m_2}{2}\right)\leq \mathcal{G}(\zeta m_1+(1-\zeta)m_2)+\mathcal{G}((1-\zeta)m_1+\zeta m_2).\label{parthhproof}
	\end{equation}
	Integrating inequality (\ref{parthhproof}) with respect to $\zeta$ over $[0,1]$, we have
	
	\begin{equation}
	\begin{aligned}	
	2\mathcal{G}\left(\frac{m_1+m_2}{2}\right)&\leq \int_{0}^{1} \mathcal{G}(\zeta m_1+(1-\zeta)m_2)d\zeta+\int_{0}^{1}\mathcal{G}((1-\zeta)m_1+\zeta m_2)d\zeta\\&=\frac{2}{m_2-m_1}\int_{m_1}^{m_2}\mathcal{G}(x)dx.
	\end{aligned}
	\end{equation}
	
	In order to prove the second part of inequality (\ref{hh}), we used Definition \ref{definitionconvexity}, for $\zeta \in [0,1]$ to arrive at
	\begin{equation*}
	\mathcal{G}(\zeta m_1+(1-\zeta)m_2)\leq \zeta\mathcal{G}(m_1)+(1-\zeta)\mathcal{G}(m_2)
	\end{equation*}
	and
	\begin{equation*}
	\mathcal{G}((1-\zeta)m_1+\zeta m_2)\leq (1-\zeta)\mathcal{G}(m_1)+\zeta\mathcal{G}(m_2).
	\end{equation*}
	When the above inequalities are added, we obtain the following
	
	\begin{equation}
	\begin{aligned}
	\mathcal{G}(\zeta m_1+&(1-\zeta)m_2)+\mathcal{G}((1-\zeta)m_1+\zeta m_2)\\&\leq \zeta\mathcal{G}(m_1)+(1-\zeta)\mathcal{G}(m_2)+(1-\zeta)\mathcal{G}(m_1)+\zeta\mathcal{G}(m_2).\label{proofhhsecondineq}
	\end{aligned}
	\end{equation}

	Integrating inequality (\ref{proofhhsecondineq}) with respect to $\zeta$ over $[0,1]$, we have
	\begin{equation*}
	\begin{aligned}
	\int_{0}^{1}	\mathcal{G}(\zeta m_1+(1-\zeta)m_2)d\zeta&+	\int_{0}^{1}\mathcal{G}((1-\zeta)m_1+\zeta m_2)d\zeta\\&\leq [\mathcal{G}(m_1)+\mathcal{G}(m_2)]\int_{0}^{1}d\zeta.
	\end{aligned}
	\end{equation*}
	Thus,
	\begin{equation*}
	\frac{2}{m_2-m_1}\int_{m_1}^{m_2}\mathcal{G}(x)dx\leq \mathcal{G}(m_1)+\mathcal{G}(m_2)
	\end{equation*}	
	completes the proof.
\end{proof}

The H-H inequality is geometrically described in \cite{niculescu2006convex}, and we have summarized it as follows:

The area under the graph of \(\mathcal{G}\) on \([m_1, m_2]\) is between the areas of two trapeziums.
While the area of the first trapezium is formed by the points of coordinates \((m_1, \mathcal{G}(m_1)),(m_2, \mathcal{G}(m_2))\) with the \(x-\)axis, that of the second trapezium is formed by the tangent to the graph of \(\mathcal{G}\) at \(\left(\frac{m_1+m_2}{2}, \mathcal{G}\left(\frac{m_1+m_2}{2}\right)\right)\) with the \(x-\)axis.

The example of H-H inequality is given as follows.

\begin{example}\cite{niculescu2004old}
	If we choose $\mathcal{G}=e^{x}$ with $x\in \mathbb{R}$, the H-H inequality yields	
	\begin{equation*}
	e^{(m_1+m_2) / 2}<\frac{e^{m_2}-e^{m_1}}{m_2-m_1}<\frac{e^{m_1}+e^{m_2}}{2},
	\end{equation*}
	for $m_1 <m_2$ in $\mathbb{R}$.

\end{example}
For more examples of H-H inequality, see \cite{khattri2010three} and \cite{dragomir2004selected}.

The importance of the H-H inequality is that each of its two sides is characterized a convex function. The necessary and sufficient condition for a continuous function \(\mathcal{G}\) to be convex on \((m_1, m_2)\) is given in the following theorem

\begin{theorem}\cite{hardy1952j}
	Let $\mathcal{G}$ be a continuous function on $(m_1, m_2)$. Then \(\mathcal{G}\) is convex iff
	\begin{equation}
	\mathcal{G}(x) \leq \frac{1}{2 z} \int_{x-z}^{x+z} \mathcal{G}(\zeta) d\zeta,\label{firstcharactrize}
	\end{equation}
	for \(m_1 \leq x-z \leq x \leq z+k \leq m_2\).
\end{theorem}

It can be shown that inequality (\ref{firstcharactrize}) is equivalent to the first part of (\ref{hh}) when \(\mathcal{G}\) is continuous on $[m_1, m_2]$  \cite{dragomir2004selected}.

The second part of inequality (\ref{hh}) can be applied as a convexity criterion in the following theorem.

\begin{theorem}\cite{robert1973convex}
	Let \(\mathcal{G}\) be continuous function on \([m_1, m_2].\) Then \(\mathcal{G}\) is convex iff
	\begin{equation*}
	\frac{1}{a_2-a_1} \int_{a_1}^{a_2} \mathcal{G}(x) d x \leq \frac{\mathcal{G}(a_1)+\mathcal{G}(a_2)}{2},
	\end{equation*}
	for all \(m_1<a_1<a_2<m_2\).
\end{theorem}

\section{H-H type inequalities for various classes of convexities}
Since different classes of convexities exist, many authors are committed to the improvements and generalizations of H-H inequality for various types of convex functions. Thus, in this section, we review some generalizations of H-H inequality involving different convex functions whose definitions are already given in Section \ref{cp1}.

Dragomir et al. \cite{dragomir1995some} established two inequalities of (\ref{hh}) which hold for classes \(Q(V)\) and 
\(P(V)\), Godunova-Levin and P-function respectively.

\begin{theorem}\cite{dragomir1995some}
	Let \(m_1, m_2 \in V\) with \(m_1<m_2\) and \(\mathcal{G} \in L_{1}[m_1, m_2].\) If \(\mathcal{G} \in Q(V),\) then
	\begin{equation}
	\mathcal{G}\left(\frac{m_1+m_2}{2}\right) \leq \frac{4}{m_2-m_1} \int_{m_1}^{m_2} \mathcal{G}(x) dx\label{H-Q(V)}
	\end{equation}
	and
	\begin{equation*}
	\frac{1}{m_2-m_1} \int_{m_1}^{m_2}\phi(x) \mathcal{G}(x) d x \leq \frac{\mathcal{G}(m_1)+\mathcal{G}(m_2)}{2}
	\end{equation*}
	hold, where \(\quad \phi(x)=\frac{(m_2-x)(x-m_1)}{(m_2-m_1)^{2}}\) and $x \in V$.
\end{theorem}
In this sense, since the constant $4$ is the best possible choice in (\ref{H-Q(V)}), it cannot be changed with any smaller constants.

\begin{theorem}\label{pfunctiontheorem}\cite{dragomir1995some}
	Let \(m_1, m_2 \in V\) with \(m_1<m_2\) and \(\mathcal{G}(x) \in L_{1}[m_1, m_2].\) If \(\mathcal{G} \in P(V),\) then
	\begin{equation}
	\mathcal{G}\left(\frac{m_1+m_2}{2}\right) \leq \frac{2}{m_2-m_1} \int_{m_1}^{m_2} \mathcal{G}(x) dx \leq 2[\mathcal{G}(m_1)+\mathcal{G}(m_2)]\label{ch2p(I)}
	\end{equation}
	holds.
\end{theorem}


For more H-H type inequalities via classes \(Q(V)\) and \(P(V)\), see \cite{pearce1999p}, \cite{barani2012hermite} and \cite{kadakal2020some}.

A variant of H-H type inequalities via $s$-convex function in second sense is proposed by Dragomir and Fitzpatrick \cite{dragomir1999hadamard}.

\begin{theorem}\label{thds}
	Suppose that $\mathcal{G}: \mathbb{R}_{+} \rightarrow \mathbb{R}_{+}$ is a $s$-convex function in the second sense, where $0<s\leq1$, $m_1, m_2 \in \mathbb{R}_{+}$ and $m_1<m_2$. If $\mathcal{G}(x) \in L_{1} [m_1, m_2]$, then
	\begin{footnotesize}
		\begin{equation}
		2^{s-1} \mathcal{G}\left(\frac{m_1+m_1}{2}\right) \leq \frac{1}{m_2-m_1} \int_{m_1}^{m_2} \mathcal{G}(x) d x \leq \frac{\mathcal{G}\left(m_1\right)+\mathcal{G}\left(m_2\right)}{s+1}\label{ch2s-convex}
		\end{equation}
	\end{footnotesize}
	holds.
\end{theorem}
The constant \(\frac{1}{s+1}\) is best possible in the second part of inequality (\ref{ch2s-convex}). We refer the reader to \cite{ozdemir2012some} and \cite{icscan2014new} for more results connected to H-H type inequalities via $s$-convex in the second sense.

Moreover, Dragomir and Fitzpatrick \cite{dragomir1999hadamard} also defined the following mapping that is closely related to (\ref{ch2s-convex}). 
\begin{theorem}
	Let \(\mathcal{G}(x):\left[m_{1}, m_{2}\right] \rightarrow \mathbb{R}\) be an $s$-convex function in the second sense on \(\left[m_{1}, m_{2}\right]\), such that \(\mathcal{G}(x) \in L_{1}[m_1, m_2]\), and $H:[0,1] \rightarrow \mathbb{R}$, then
	\begin{footnotesize}
		\begin{equation*}
		H(\zeta)=\frac{1}{m_{2}-m_{1}} \int_{m_{1}}^{m_{2}} \mathcal{G}\left(\zeta x+(1-\zeta) \frac{m_{1}+m_{2}}{2}\right) dx
		\end{equation*}
	\end{footnotesize}
	holds, for $\zeta\in [0,1]$.
\end{theorem}
The properties of the mapping $H$ are given as follows:

\begin{itemize}
	\item [i.] $H \in {K}_{s}^{2}$ on \([0,1]\),
	\item[ii.] $H \geq 2^{s-1} \mathcal{G}(\frac{m_1+m_2}{2})$.
\end{itemize}

These properties are the generalization of some results from \cite{dragomir1991mapping}. Also, for more properties of mappings associated with H-H inequality, see \cite{dragomir1992two}, \cite{dragomir1992some}, \cite{dragomir1993ds}, \cite{dragomir2002new}, \cite{dragomir2003further} and \cite{kilicman2015some}.


Another new H-H type inequalities for the preinvex function is given by Noor \cite{noor2007hermite} as follows.

\begin{theorem}\label{ch2thepreinvexity}
	Let $\mathcal{G}:V=[m_1, m_1+\eta (m_2,m_1)]\rightarrow (0, \infty)$ be a preinvex function on $V^{\circ}$ with \(m_1<m_1+\eta (m_2,m_1)\), $m_1,m_2 \in V^{\circ}$  and \(\mathcal{G} \in L_{1}[m_1, m_2],\) then
	\begin{footnotesize}
		\begin{equation}
		\mathcal{G}\left(\frac{2m_1+\eta (m_2,m_1)}{2}\right) \leq \frac{1}{\eta (m_2,m_1)} \int_{m_1}^{m_1+\eta (m_2,m_1)} \mathcal{G}(x) dx \leq\frac{\mathcal{G}(m_1)+\mathcal{G}(m_2)}{2}.
		\label{ch2ineqpreinvexity}
		\end{equation}
	\end{footnotesize}
\end{theorem}

\begin{remark}
	In Theorem \ref{ch2thepreinvexity}, if we take $\eta (m_2,m_1)=m_2-m_1$, then inequality (\ref{ch2ineqpreinvexity}) reduces to (\ref{hh}).	
\end{remark}

\section{H-H type inequalities for differentiable functions}\label{ch2sec2.3}

An interesting problem in (\ref{hh}) that attracts many researchers is the determination of two bounds of quantities (\ref{Interestingone}) and (\ref{Interestingtwo}) given as follows:

\begin{equation}
\bigg|\frac{1}{m_2-m_1} \int_{m_1}^{m_2} \mathcal{G}(x) d x-\mathcal{G}\bigg(\frac{m_1+m_2}{2}\bigg)\bigg|\label{Interestingone}
\end{equation}
and
\begin{equation}
\bigg|\frac{\mathcal{G}(m_1)+\mathcal{G}(m_2)}{2}-\frac{1}{m_2-m_1} \int_{m_1}^{m_2} \mathcal{G}(x) d x\bigg|.\label{Interestingtwo}
\end{equation}

While (\ref{Interestingone}) estimates the difference between the left and the middle parts of (\ref{hh}), the quantity (\ref{Interestingtwo}) estimates the difference between the middle and the right parts of (\ref{hh}). The quantity (\ref{Interestingone}) is called the mid-point type inequality. Meanwhile, the 
quantity (\ref{Interestingtwo}) is named as the trapezoid type inequality. Recently, different integral inequalities were obtained through differentiable convexity. The following result, given by Dragomir and Agarwal \cite{dragomir1998two}, can be used to estimate a new bound on (\ref{Interestingtwo}).

\begin{lemma}\label{ch2lamma1}
	Let \(\mathcal{G}: V \subseteq \mathbb{R} \rightarrow \mathbb{R}\) be a differentiable function on \(V^{\circ}, m_1, m_2 \in V^{\circ}\) with \(m_1<m_2\) and \(\mathcal{G}^{\prime} \in L_{1}[m_1, m_2],\) then the following identity holds:
	\begin{footnotesize}
		\begin{equation*}
		\frac{\mathcal{G}(m_1)+\mathcal{G}(m_2)}{2}-\frac{1}{m_2-m_1} \int_{m_1}^{m_2} \mathcal{G}(x) dx=\frac{m_2-m_1}{2} \int_{0}^{1}(1-2 \zeta) \mathcal{G}^{\prime}(\zeta m_1+(1-\zeta) m_2) d\zeta.
		\end{equation*}	
	\end{footnotesize}
\end{lemma}

Therefore, using Lemma \ref{ch2lamma1} the following theorems connected with the second part of (\ref{hh}) for differentiable convex functions hold.

\begin{theorem}\label{ch2theorem1}
	Let \(\mathcal{G}: V \subseteq \mathbb{R} \rightarrow \mathbb{R}\) be a differentiable function on \(V^{\circ}, m_1, m_2 \in V^{\circ}\) with \(m_1<m_2\) and \(\mathcal{G}^{\prime} \in L_{1}[m_1, m_2],\) if \(\left|\mathcal{G}^{\prime}\right|\) is convex on \([m_1, m_2],\) then
	\begin{footnotesize}
		\begin{equation}
		\left|\frac{\mathcal{G}(m_1)+\mathcal{G}(m_2)}{2}-\frac{1}{m_2-m_1} \int_{m_1}^{m_2} \mathcal{G}(x) dx\right| \leq \frac{m_2-m_1}{8}\left[\left|\mathcal{G}^{\prime}(m_1)\right|+\left|\mathcal{G}^{\prime}(m_2)\right|\right].\label{ch2ineq1sec2.3}
		\end{equation}
	\end{footnotesize}
\end{theorem}

\begin{theorem}\label{ch2theorem2}
	Let \(\mathcal{G}: V \subseteq \mathbb{R} \rightarrow \mathbb{R}\) be a differentiable function on \(V^{\circ}, m_1, m_2 \in V^{\circ}\) with \(m_1<m_2\) and \(\mathcal{G}^{\prime} \in L_{1}[m_1, m_2].\) If \(\left|\mathcal{G}^{\prime}\right|^{q}\) is convex on \([m_1, m_2]\) for \(q>1\) with
	\(q(p-1)=p,\) then
	\begin{footnotesize}
		\begin{equation}
		\left|\frac{\mathcal{G}(m_1)+\mathcal{G}(m_2)}{2}-\frac{1}{m_2-m_1} \int_{m_1}^{m_2} \mathcal{G}(x) d x\right| \leq \frac{m_2-m_1}{2}\left(\frac{1}{p+1}\right)^{\frac{1}{p}}\left[\frac{\left|\mathcal{G}^{\prime}(m_1)\right|^{q}+\left|\mathcal{G}^{\prime}(m_2)\right|^{q}}{2}\right]^{\frac{1}{q}}.
		\end{equation}
	\end{footnotesize}
\end{theorem}

The improvement and simplification of the aforementioned result presented in Theorem \ref{ch2theorem2} is provided by Pearce and Pe\u{c}ari\'{c} \cite{pearce2000inequalities}.

\begin{theorem}\label{ch2theorem3}
	Let \(\mathcal{G}: V \subseteq \mathbb{R} \rightarrow \mathbb{R}\) be a differentiable function on $V^{\circ}$, with $m_1, m_2 \in V^{\circ}$,  $m_1<m_2$ and \(\mathcal{G}^{\prime} \in L_{1}[m_1, m_2].\) If \(\left|\mathcal{G}^{\prime}\right|^{q}\), for \(q>1\) where
	\(q(p-1)=p,\) is convex on \([m_1, m_2],\) then
	\begin{footnotesize}
		\begin{equation}
		\left|\frac{\mathcal{G}(m_1)+\mathcal{G}(m_2)}{2}-\frac{1}{m_2-m_1} \int_{m_1}^{m_2} \mathcal{G}(x) d x\right| \leq \frac{m_2-m_1}{4}\left[\frac{\left|\mathcal{G}^{\prime}(m_1)\right|^{q}+\left|\mathcal{G}^{\prime}(m_2)\right|^{q}}{2}\right]^{\frac{1}{q}}.
		\end{equation}
	\end{footnotesize}
\end{theorem}

\begin{remark}
	Choosing \(q=1\) reduces Theorem \ref{ch2theorem3} to Theorem \ref{ch2theorem1}. In Theorem \ref{ch2theorem3}, taking $q=\frac{p}{p-1}$ improves the constant given in Theorem \ref{ch2theorem2} since 
	\(\frac{1}{4}<\frac{1}{2(p+1)^{\frac{1}{p}}}\), where $p>1$.	
\end{remark}

Kirmaci \cite{kirmaci2004inequalities} proved the following results that give the bounds on (\ref{Interestingone}) by using the assumptions of convexity.

\begin{lemma}\label{ch2lemmamidhalf}
	Let $\mathcal{G} : V \subseteq \mathbb{R} \rightarrow \mathbb{R}$ be a differentiable mapping on $V^{\circ}$, $m_1, m_2 \in V^{\circ}$ with $m_1<m_2$. If \(\mathcal{G}^{\prime} \in L_{1}[m_1, m_2],\) then we have
	\begin{footnotesize}
		\begin{equation*}
		\frac{1}{m_2-m_1} \int_{m_1}^{m_2} \mathcal{G}(x) dx-\mathcal{G}\left(\frac{m_1+m_2}{2}\right)=(m_2-m_1) \int_{0}^{1} Q(\zeta) \mathcal{G}^{\prime}(\zeta m_1+(1-\zeta) m_2) d\zeta,
		\end{equation*}
	\end{footnotesize}
	where
	\begin{footnotesize}
		\begin{equation*}
		Q(\zeta)=\left\{\begin{array}{ll}{\zeta,} & {\zeta \in\left[0, \frac{1}{2}\right)} \\ {\zeta-1,} & {\zeta \in\left[\frac{1}{2}, 1\right].}\end{array}\right.
		\end{equation*}
	\end{footnotesize}
\end{lemma}

\begin{theorem}\label{thmidclassch2}
	Let $\mathcal{G} : V\subseteq \mathbb{R}$ be a differentiable mapping on $V^{\circ}$, $m_1, m_2 \in V^{\circ}$ with $m_1<m_2$. If $\left|\mathcal{G}^{\prime}\right|$ is convex on $[m_1, m_2]$, then we have
	\begin{footnotesize}
		\begin{flalign}
		{\quad\left|\frac{1}{m_2-m_1} \int_{m_1}^{m_2} \mathcal{G}(x) dx-\mathcal{G}\left(\frac{m_1+m_2}{2}\right)\right| \leqslant \frac{m_2-m_1}{8}\left(\left|\mathcal{G}^{\prime}(m_1)\right|+\left|\mathcal{G}^{\prime}(m_2)\right|\right)}.\label{ineqmidclassch2}
		\end{flalign}
	\end{footnotesize}	
\end{theorem}

Some new inequalities for twice differentiable functions connected to inequality (\ref{hh}) are given by Dragomir and Pearce \cite{dragomir2004selected} through the following lemma.

\begin{lemma}\label{dragomirsecondclassical}
	Let \(\mathcal{G}: V \subseteq \mathbb{R} \rightarrow \mathbb{R}\) be a twice differentiable function on \(V^{\circ}, m_1, m_2 \in V^{\circ}\) with \(m_1<m_2\) and \(\mathcal{G}^{\prime\prime } \in L_{1}[m_1, m_2],\) then the following
	\begin{footnotesize}
		\begin{equation*}
		\frac{\mathcal{G}(m_1)+\mathcal{G}(m_2)}{2}-\frac{1}{m_2-m_1} \int_{m_1}^{m_2} \mathcal{G}(x) d x=\frac{(m_2-m_1)^2}{2} \int_{0}^{1}\zeta(1-\zeta) \mathcal{G}^{\prime\prime }(\zeta m_1+(1-\zeta) m_2) d\zeta
		\end{equation*}	
	\end{footnotesize}
	holds.
\end{lemma}

Kirmaci et al. \cite{kirmaci2007hadamard} studied new inequality of H-H type for differentiable mappings involving $s$-convexity.  

\begin{theorem} Let \(\mathcal{G}: V \subset[0, \infty) \rightarrow \mathbb{R}\) be a differentiable
	mapping on \(V^{\circ}\) such that \(\mathcal{G}^{\prime} \in L_{1}[m_1, m_2],\) where \(m_1, m_2 \in V\) with \(m_1<m_2 .\) If \(\left|\mathcal{G}^{\prime}\right|^{q}\) is $s$-convex on
	\([m_1, m_2],\) where $q \geq 1$ and \(s \in(0,1],\) we have
	\begin{footnotesize}
		\begin{equation}
		\begin{aligned} \bigg| \frac{\mathcal{G}(m_1)+\mathcal{G}(m_2)}{2} &-\frac{1}{m_2-m_1} \int_{m_1}^{m_2} \mathcal{G}(x) dx \bigg| \\ \leq & \frac{m_2-m_1}{2}\left(\frac{1}{2}\right)^{\frac{q-1}{q}}\left[\frac{s+\left(\frac{1}{2}\right)^{s}}{(s+1)(s+2)}\right]^{\frac{1}{q}}\left(\left|\mathcal{G}^{\prime}(m_1)\right|^{q}+\left|\mathcal{G}^{\prime}(m_2)\right|^{q}\right)^{\frac{1}{q}}. \end{aligned}
		\end{equation}
	\end{footnotesize}
\end{theorem}

Barani et al. \cite{baranietal2012hermite} generalized Lemma \ref{ch2lamma1} to estimate the trapezoid type inequalities connected with (\ref{hh}) for preinvex function.

\begin{lemma}\label{ch2lemma}
	Suppose that $\mathcal{G}:V=[m_1,m_1+\eta(m_2,m_1)]\to (0,\infty)$ is a differentiable function, where $m_1,m_1+\eta(m_2,m_1)\in V$ with $m_1<m_1+\eta(m_2,m_1)$. If $\mathcal{G}^{\prime}\in L_{1}[m_1,m_1+\eta(m_2,m_1)]$, we have
	\begin{footnotesize}
		\begin{flalign*}
		\frac{1}{\eta(m_2,m_1)}\int_{m_1}^{m_1+\eta(m_2,m_1)}\mathcal{G}(x)dx-\frac{\mathcal{G}(m_1)+\mathcal{G}(m_1+\eta(m_2,m_1))}{2}\\=\frac{\eta(m_2,m_1)}{2}\bigg[\int_{0}^{1}(1-2\zeta)\mathcal{G}^{\prime}(m_1+\zeta\eta(m_2,m_1))d\zeta\bigg].\label{lem1}
		\end{flalign*}	
	\end{footnotesize}
\end{lemma}

Recently, presumably new H-H type inequalities were established by Mehrez and Agarwal \cite{mehrez2019new}, whose finding is reported in the next theorem.

\begin{theorem}\label{con2}
	Suppose that $\mathcal{G}:V\subseteq\mathbb{R}\to\mathbb{R}$ is a differentiable mapping on $V^{\circ}$, $m_1,m_2\in V^{\circ}$ with $m_1<m_2$. Let the derivative of $\mathcal{G}$ be $\mathcal{G}^{\prime}:[\frac{3m_1-m_2}{2},\frac{3m_2-m_1}{2}]\to \mathbb{R}$, a continuous function on $[\frac{3m_1-m_2}{2},\frac{3m_2-m_1}{2}]$. Let $q\geq1$, if $|\mathcal{G}^{\prime}|$ is convex on $[\frac{3m_1-m_2}{2},\frac{3m_2-m_1}{2}]$, then the following
	\begin{footnotesize}
		\begin{equation}
		\begin{aligned}
		\bigg|\frac{1}{m_2-m_1}\int_{m_1}^{m_2}\mathcal{G}(x)dx&-\mathcal{G}\bigg(\frac{m_1+m_2}{2}\bigg)\bigg|\\&\leq\frac{m_2-m_1}{8}\bigg
		(\bigg|\mathcal{G}^{\prime}\bigg(\frac{3m_1-m_2}{2}\bigg)\bigg|^{q}+\bigg|\mathcal{G}^{\prime}\bigg(\frac{3m_1-m_2}{2}\bigg)\bigg|^{q}\bigg)^{\frac{1}{q}}\label{gh}
		\end{aligned}
		\end{equation}
	\end{footnotesize}
	holds.
\end{theorem}

Almutairi and Kili\c{c}man \cite{almutairi2019new} extended Theorem \ref{con2} to $s$-convex in the second sense as follows.

\begin{theorem}
	Suppose that $\mathcal{G}:V\subseteq\mathbb{R}_{+}\to\mathbb{R}$ is a differentiable mapping on $V^{\circ}$, $m_1,m_2\in V^{\circ}$ with $m_1<m_2$. Let the derivative of $\mathcal{G}$ be $\mathcal{G}^{\prime}:[\frac{3m_1-m_2}{2},\frac{3m_2-m_1}{2}]\to \mathbb{R}$, a continuous function on $[\frac{3m_1-m_2}{2},\frac{3m_2-m_1}{2}]$. Let $q\geq1$, if $|\mathcal{G}^{\prime}|^{q}$ is an $s$-convex function on $[\frac{3m_1-m_2}{2},\frac{3m_2-m_1}{2}]$, for some fixed $s\in(0,1)$, then we have the following:
	\begin{footnotesize}	
		\begin{flalign}
		\bigg|&\frac{1}{m_2-m_1}\int_{m_1}^{m_2}\mathcal{G}(x)dx-\mathcal{G}\bigg(\frac{m_1+m_2}{2}\bigg)\bigg|\nonumber\\&\leq (m_2-m_1) \bigg(\frac{1}{8}\bigg)^{\frac{q-1}{q}} \bigg[\frac{2-2^{-s}}{2(s+1)(s+2)}\bigg]^\frac{1}{q}\bigg[\bigg|\mathcal{G}^{\prime}\bigg(\frac{3m_1-m_2}{2}\bigg)\bigg|^{q}+\bigg|\mathcal{G}^{\prime}\bigg(\frac{3m_2-m_1}{2}\bigg)\bigg|^{q}\bigg]^{\frac{1}{q}}.\label{eqs1}
		\end{flalign}
	\end{footnotesize}	
\end{theorem}

\section{Generalized H-H type inequalities involving different fractional integrals}

This section presents some results on the generalization of inequalities introduced in Section \ref{ch2sec2.3}. Therefore, many generalizations of H-H type inequalities established using fractional integrals for different classes of convexities are discussed here since they can be frequently used in other parts of the thesis. For example, the work of Sarikaya et al. \cite{sarikaya2013hermite} was the first to present inequalities of H-H type involving Riemann-Liouville fractional integrals. This is given below.
\begin{theorem}\label{smth}
	Suppose that $\mathcal{G} :[m_1, m_2] \rightarrow \mathbb{R}$  is a non-negative function with $0 \leq m_1<m_2$  and $\mathcal{G} \in L_{1}[m_1, m_2]$. If $\mathcal{G}$ is a convex function on $[m_1, m_2]$, we have
	\begin{footnotesize}
		\begin{equation}
		\mathcal{G}\left(\frac{m_1+m_2}{2}\right) \leq \frac{\Gamma(\lambda+1)}{2(m_2-m_1)^{\lambda}}\left[J_{m_1^{+}}^{\lambda}\mathcal{G}(m_2)+J_{m_2^{-}}^{\lambda}\mathcal{G}(m_1)\right] \leq \frac{\mathcal{G}(m_1)+\mathcal{G}(m_2)}{2},\label{smin}
		\end{equation}
	\end{footnotesize}
	where $\lambda>0$.
\end{theorem}

\begin{remark}
	In Theorem \ref{smth}, choosing $\lambda=1$ reduces inequality (\ref{smin}) to (\ref{hh}).
\end{remark}

Moreover, Sarikaya et al. \cite{sarikaya2013hermite} presented the following fractional integral identity. 

\begin{lemma}\label{ch2sec24lemmaone}
	Let  $\mathcal{G} :[m_1, m_2] \rightarrow \mathbb{R}$  be a differentiable function on $(m_1, m_2)$  with  $m_1<m_2$. If $\mathcal{G}^{\prime} \in L_{1}[m_1, m_2]$, then we have
	\begin{footnotesize}
		\begin{equation*}
		\begin{aligned}
		\frac{\mathcal{G}(m_1)+\mathcal{G}(m_2)}{2}-&\frac{\Gamma(\lambda+1)}{2(m_2-m_1)^{\lambda}}\left[J_{m_1^+}^{\lambda} \mathcal{G}(m_2)+J_{m_2^-}^{\lambda} \mathcal{G}(m_1)\right]\\=&\frac{m_2-m_1}{2} \int_{0}^{1}\left[(1-\zeta)^{\lambda}-\zeta^{\lambda}\right] \mathcal{G}^{\prime}(\zeta m_1+(1-\zeta) m_2) d \zeta.
		\end{aligned}
		\end{equation*}
	\end{footnotesize}
\end{lemma}

The identity presented in the above lemma was also used by Sarikaya when determining trapezoid type inequalities connected with (\ref{hh}) for Riemann-Liouville fractional integrals.
\begin{theorem}\label{sfth}
	Let  $\mathcal{G} :[m_1, m_2] \rightarrow \mathbb{R}$  be a differentiable function on $(m_1, m_2)$  with  $m_1<m_2$ and $\mathcal{G}^{\prime} \in L_{1}[m_1, m_2]$. If \(|\mathcal{G}^{\prime} |\) is convex on \([m_1, m_2]\), then we have
	\begin{footnotesize}
		\begin{equation}
		\begin{aligned}
		\bigg|\frac{\mathcal{G}(m_1)+\mathcal{G}(m_2)}{2}-&\frac{\Gamma(\lambda+1)}{2(m_2-m_1)^{\lambda}}\left[J_{m_1^+}^{\lambda} \mathcal{G}(m_2)+J_{m_2^-}^{\lambda} \mathcal{G}(m_1)\right]\bigg|\\ \leq& \frac{m_2-m_1}{2(\lambda+1)}\left(1-\frac{1}{2^{\lambda}}\right)\left[\mathcal{G}^{\prime}(m_1)+\mathcal{G}^{\prime}(m_2)\right].\label{sfin}
		\end{aligned}
		\end{equation}
	\end{footnotesize}
\end{theorem}

\begin{remark}
	Taking $\lambda=1$ in Theorem \ref{sfth} reduces inequality (\ref{sfin})
	to inequality (\ref{ch2ineq1sec2.3}) of Theorem \ref{ch2theorem1}.
\end{remark}

Zhu et al. \cite{zhu2012fractional} studied a new fractional integral identity for differentiable convex mappings. The results are presented below.
\begin{lemma}\label{zhulemma}
	Let \(\mathcal{G}:[m_1, m_2] \rightarrow \mathbb{R}\) be a differentiable mapping on \((m_1, m_2)\) with \(m_1<m_2.\) If \(\mathcal{G}^{\prime} \in\)
	\(L_{1}[m_1, m_2],\) then the equality for fractional integrals holds as follows:
	\begin{footnotesize}
		\begin{equation}
		\begin{aligned}
		&\frac{\Gamma(\lambda+1)}{2(m_2-m_1)^{\lambda}}\left[J_{m_1^+}^{\lambda} \mathcal{G}(m_2)+J_{m_2^-}^{\lambda} \mathcal{G}(m_1)\right]-\mathcal{G}\left(\frac{m_1+m_2}{2}\right)\\&= \frac{m_2-m_1}{2}\left[\int_{0}^{1} P  \mathcal{G}^{\prime}(\zeta m_1+(1-\zeta) m_2) d \zeta - \int_{0}^{1}\left[(1-\zeta)^{\lambda}-\zeta^{\lambda}\right] \mathcal{G}^{\prime}(\zeta m_1+(1-\zeta) m_2) d \zeta\right],\label{identityzhu} \end{aligned}
		\end{equation}
	\end{footnotesize}
	where
	\begin{equation*}
	P=\left\{\begin{array}{ll}{1,} & {0 \leq \zeta<\frac{1}{2}} \\ {-1,} & {\frac{1}{2} \leq \zeta<1.}\end{array}\right.
	\end{equation*}
\end{lemma}

Using the above identity, the following result estimates a midpoint type inequalities related to (\ref{hh}), which involves Riemann-Liouville fractional integrals. 

\begin{theorem}\label{zhutheorem}
	Let \(\mathcal{G}:[m_1, m_2] \rightarrow \mathbb{R}\) be a differentiable mapping on \((m_1, m_2)\) with \(m_1<m_2.\) If \(\left|\mathcal{G}^{\prime}\right|\) is convex on \([m_1, m_2],\)
	then the following inequality holds:
	\begin{footnotesize}
		\begin{equation}
		\begin{aligned}&
		\left|\frac{\Gamma(\lambda+1)}{2(m_2-m_1)^{\lambda}}\left[J_{m_1^+}^{\lambda} \mathcal{G}(m_2)+J_{m_2^-}^{\lambda} \mathcal{G}(m_1)\right]-\mathcal{G}\left(\frac{m_1+m_2}{2}\right)\right|\\\leq& \frac{m_2-m_1}{4(\lambda+1)}\left(\lambda+3-\frac{1}{2^{\lambda-1}}\right)\left[\left|\mathcal{G}^{\prime}(m_1)\right|+\left|\mathcal{G}^{\prime}(m_2)\right|\right].\label{ineqzho}
		\end{aligned}
		\end{equation}
	\end{footnotesize}
\end{theorem}

Almutair and Kili\c{c}man \cite{almutairi2020new} extended Lemma \ref{zhulemma} and Theorem \ref{zhutheorem} for Katugampola fractional integrals as follows.

\begin{lemma}\label{lemmain}
	Let \(\mathcal{G}:[m_1^{\rho}, m_2^{\rho}] \rightarrow \mathbb{R}\) be a differentiable mapping on \((m_1^{\rho}, m_2^{\rho})\), where \(m_1<m_2.\) The following equality holds if the fractional integrals exist,
	
	\begin{footnotesize}	
		\begin{equation}
		\begin{aligned} &\frac{\lambda\rho^{\lambda}\Gamma(\lambda+1)}{2(m_2^{\rho}-m_1^{\rho})^{\lambda}}[^{\rho} I_{m_1^+}^{\lambda}\mathcal{G}(m_2^\rho )+^{\rho} I_{m_2^-}^{\lambda}\mathcal{G}(m^\rho )]-\mathcal{G}\left(\frac{m_1^{\rho}+m_2^{\rho}}{2}\right) \\=& \frac{m_2^{\rho}-m_1^{\rho}}{2}\left[\int_{0}^{1} M \mathcal{G}^{\prime}(\zeta^{\rho} m_1^{\rho}+(1-\zeta^{\rho}) m_2^{\rho}) d \zeta-\int_{0}^{1}\left[(1-\zeta^\rho)^{\lambda}-\zeta^{\rho \lambda}\right]\zeta^{\rho -1} \mathcal{G}^{\prime}(\zeta^{\rho} m_1^{\rho}+(1-\zeta^{\rho}) m_2^{\rho}) d \zeta\right],\label{identi} \end{aligned}
		\end{equation}
	\end{footnotesize}
	where
	\begin{equation*}
	M=\left\{\begin{array}{ll}{\zeta^{\rho -1},} & {0 \leq \zeta<\frac{1}{^{\rho}\sqrt{2}}} \\ {-\zeta^{\rho -1},} & {\frac{1}{^{\rho}\sqrt{2}} \leq \zeta<1}.\end{array}\right.
	\end{equation*}
\end{lemma}

\begin{remark}
	If \(\rho=1,\) then the identity (\ref{identi}) in Lemma \ref{lemmain} reduces to identity (\ref{identityzhu}) in Lemma \ref{zhulemma}.
\end{remark}

\begin{theorem}\label{theorem}
	Let \(\mathcal{G}:\left[m_1^{\rho}, m_2^{\rho}\right] \rightarrow \mathbb{R}\) be a differentiable mapping on $(m_1^{\rho}, m_2^{\rho})$ with \(0 \leq m_1<m_2 .\) If \(\left|\mathcal{G}^{\prime}\right|\) is convex
	on \(\left[m_1^{\rho}, m_2^{\rho}\right],\) then the following inequality
	\begin{footnotesize}	
		\begin{equation}
		\begin{aligned}&
		\left|\frac{\lambda \rho^{\lambda} \Gamma(\lambda+1)}{2\left(m_2^{\rho}-m_1^{\rho}\right)^{\lambda}}\left[^{\rho} I_{m_1^+}^{\lambda} \mathcal{G}\left(m_2^{\rho}\right)+^{\rho} I_{m_2^-}^{\lambda} \mathcal{G}\left(m_1^{\rho}\right)\right]-\mathcal{G}\bigg(\frac{m_1^{\rho}+m_2^{\rho}}{2}\bigg)\right|\\ \leq& \frac{m_2^{\rho}-m_1^{\rho}}{4 \rho(\lambda +1)}\bigg[3+\lambda -\frac{1}{2^{\lambda -1}}\bigg]\left[\left|\mathcal{G}^{\prime}\left(m_1^{\rho}\right)\right|+\left|\mathcal{G}^{\prime}\left(m_2^{\rho}\right)\right|\right]\label{mainineq}
		\end{aligned}
		\end{equation}
		holds.
	\end{footnotesize}
\end{theorem}

\begin{remark}	
	Consider inequality (\ref{mainineq}) of the Theorem \ref{theorem}, we have
	\begin{itemize}	
		\item[i.] Choosing $\rho=1$  reduces inequality (\ref{mainineq}) to inequality (\ref{ineqzho}) of Theorem \ref	{zhutheorem}.	
		\item[ii.] Taking $\rho=1$ and $\lambda=1$ reduces inequality (\ref{mainineq}) to inequality (16) in \cite{zhu2012fractional}, which is given as follows
		\begin{footnotesize}
			\begin{equation*}
			\left|\frac{1}{m_2-m_1} \int_{m_1}^{m_2} \mathcal{G}(x) d x-\mathcal{G}\left(\frac{m_1+m_2}{2}\right)\right| \leq \frac{3(m_2-m_1)}{8}\left(\left|\mathcal{G}^{\prime}(m_1)\right|+\left|\mathcal{G}^{\prime}(m_2)\right|\right).
			\end{equation*}
		\end{footnotesize}
	\end{itemize}
\end{remark}

Meanwhile, Wang et al. \cite{wang2013hadamard} extended Lemma \ref{ch2sec24lemmaone} to include a twice-differentiable mapping. 

\begin{lemma}\label{lemsed}
	Let $\mathcal{G}:[m_1, m_2] \rightarrow \mathbb{R}$  be a twice-differentiable function on $(m_1, m_2)$  with $m_1<m_2$. If $\mathcal{G}^{\prime \prime} \in L_{1}[m_1, m_2]$, then
	\begin{footnotesize}	
		\begin{equation*}
		\begin{aligned} & \frac{\mathcal{G}(m_1)+\mathcal{G}(m_2)}{2}-\frac{\Gamma(\lambda+1)}{2(m_2-m_1)^{\lambda}}\left[J_{m_1^+}^{\lambda} \mathcal{G}(m_2)+J_{m_2^-}^{\lambda} \mathcal{G}(m_1)\right] \\=& \frac{(m_2-m_1)^{2}}{2} \int_{0}^{1} \frac{1-(1-\zeta)^{\lambda+1}-\zeta^{\lambda+1}}{\lambda+1} \mathcal{G}^{\prime \prime}(\zeta m_1+(1-\zeta) m_2) d \zeta \end{aligned}
		\end{equation*}
	\end{footnotesize}
	holds.
\end{lemma}

Set et al. \cite{set2014hermite} generalized Theorem \ref{thds} for fractional integrals, and the result is given as follows.

\begin{theorem}\label{ths_Rch2}
	Suppose that $\mathcal{G} :[m_1, m_2] \rightarrow \mathbb{R}$  is a non-negative function with $0 \leq m_1<m_2$  and $\mathcal{G} \in L_{1}[m_1, m_2]$. If $\mathcal{G}$ is $s$-convex function in the second sense on $[m_1, m_2]$, we have
	\begin{footnotesize}
		\begin{equation}
		\begin{aligned}
		2^{s-1}	\mathcal{G}\left(\frac{m_1+m_2}{2}\right) &\leq \frac{\Gamma(\lambda+1)}{2(m_2-m_1)^{\lambda}}\left[J_{m_1^{+}}^{\lambda}\mathcal{G}(m_2)+J_{m_2^{-}}^{\lambda}\mathcal{G}(m_1)\right] \\&\leq \bigg[\frac{\lambda}{(\lambda+s)}+\lambda\beta(\lambda,s+1)\bigg] \frac{\mathcal{G}(m_1)+\mathcal{G}(m_2)}{2},\label{Ins_Rch2}
		\end{aligned}
		\end{equation}
	\end{footnotesize}
	where $\lambda>0$ and $0<s<1$.
\end{theorem}

The H-H inequality for Hadamard fractional integrals that was established by Wang et al. \cite{wang2013refinements} also received the attention of many researchers. This refinement is given as follows.

\begin{theorem}
	Suppose that \(\mathcal{G}:[m_1, m_2] \rightarrow \mathbb{R}\) is a non-negative function with \(0<m_1<m_2\) and \(\mathcal{G} \in L_{1}[m_1, m_1].\)
	If \(\mathcal{G}\) is a non-decreasing convex function on \([m_1, m_2],\) then the following inequality 
	\begin{equation*}
	\mathcal{G}(\sqrt{m_1 m_2}) \leq \frac{\Gamma(\lambda+1)}{2(\ln m_2-\ln m_1)^{\lambda}}\left[ H_{m_1^{+}}^{\lambda} \mathcal{G}(m_2)+ H_{m_2^{-}}^{\lambda} \mathcal{G}(m_1)\right] \leq \mathcal{G}(m_2)
	\end{equation*}
	holds.
\end{theorem}

Mo et al. \cite{mo2014generalized} provided the generalized H-H type inequalities involving local fractional integrals for generalized convex function on fractal sets as follows.

\begin{theorem}
	Let \(\mathcal{G}(x) \in I_{x}^{\alpha}[m_1, m_2]\)
	be a generalized convex function on \([m_1, m_2]\) with \(m_1<m_2.\) Then
	\begin{equation}
	\mathcal{G}\left(\frac{m_1+m_2}{2}\right) \leq \frac{\Gamma(1+\alpha)}{(m_2-m_1)^{\alpha}} m_1 I_{m_2}^{\alpha} \mathcal{G}(x) \leq \frac{\mathcal{G}(m_1)+\mathcal{G}(m_2)}{2^{\alpha}}\label{localhh}
	\end{equation}
	holds.
\end{theorem}

\begin{remark}
	Choosing $\alpha=1$ in inequality (\ref{localhh}), we get inequality (\ref{hh}).
\end{remark}

Furthermore, the H-H type inequalities for the generalized $s$-convex function in the second sense on fractal sets were proposed by Mo and Sui \cite{mo2017hermite}.

\begin{theorem}
	Suppose that \(\mathcal{G}: \mathbb{R}_{+} \rightarrow \mathbb{R}^{\alpha}\) is a generalized s-convex function in the second sense for \(0<s<1\) and \(m_1, m_2 \in[0, \infty)\) with \(m_1<m_2.\) Then, for \(\mathcal{G} \in C_{\alpha}[m_1, m_2],\) the following inequality
	\begin{footnotesize}
		\begin{equation}
		\frac{2^{(s-1) \alpha}}{\Gamma(1+\alpha)}\mathcal{G}\left(\frac{m_1+m_2}{2}\right) \leq \frac{m_1 I_{m_2}^{\alpha} \mathcal{G}(x)}{(m_2-m_1)^{\alpha}} \leq \frac{\Gamma(1+s \alpha)}{\Gamma(1+(s+1) \alpha)}(\mathcal{G}(m_1)+\mathcal{G}(m_2))\label{s-convexlocal}
		\end{equation}
	\end{footnotesize}
	holds.
\end{theorem}
\begin{remark}
	When taking $\alpha=1$ in (\ref{s-convexlocal}), we obtained (\ref{ch2s-convex}).
\end{remark}

For more results on the generalizations of H-H type inequalities involving fractal sets via fractional integrals, one should consult the following references \cite{vivas2016new,luo2020fejer}.

The result in Theorem \ref{ch2sec24endtheorem}, involving Katugampola fractional integrals, is the generalization of the result earlier presented in Theorem \ref{smth}.

\begin{theorem}\cite{chen2017hermite}\label{ch2sec24endtheorem}
	Let \(\lambda>0\) and \(\rho>0.\) Let \(\mathcal{G}:\left[m_1^{\rho}, m_2^{\rho}\right] \rightarrow \mathbb{R}\) be a positive function with \(0 \leq m_1<m_2\) and \(\mathcal{G} \in X_{c}^{p}\left(m_1^{\rho}, m_2^{\rho}\right).\) If \(	\mathcal{G}\) is also a convex function on \([m_1, m_2],\) then the following inequality
	\begin{footnotesize}
		\begin{equation}
		\mathcal{G}\left(\frac{m_1^{\rho}+m_2^{\rho}}{2}\right) \leq \frac{\rho^{\lambda} \Gamma(\lambda+1)}{2\left(m_2^{\rho}-m_1^{\rho}\right)^{\lambda}}\left[^{\rho} I_{m_1^+}^{\lambda} \mathcal{G}\left(m_2^{\rho}\right)+^{\rho} I_{m_2^-}^{\lambda} \mathcal{G}\left(m_1^{\rho}\right)\right] \leq \frac{\mathcal{G}\left(m_1^{\rho}\right)+\mathcal{G}\left(m_2^{\rho}\right)}{2}\label{ch2sec24cheninequality}
		\end{equation}
	\end{footnotesize}
	holds, where the fractional integrals are considered for the function \(\mathcal{G}\left(x^{\rho}\right)\) and evaluated at $m_1$ and $m_2$, respectively.
\end{theorem}

The estimate of the difference between the right term and the middle term of inequality (\ref{ch2sec24cheninequality}) is obtained using the following lemma.

\begin{lemma}\label{lem}
	Suppose that $\mathcal{G}: \left[m_1^{\rho}, m_2^{\rho}\right] \subset \mathbb{R}_{+} \rightarrow \mathbb{R}$ is a differentiable mapping on $\left(m_1^{\rho}, m_2^{\rho}\right)$, where $0 \leq m_1<m_2$ and $\lambda,\rho>0$. If the fractional integrals exist, we have
	\begin{footnotesize}
		\begin{equation}
		\begin{aligned}
		&	\frac{\mathcal{G}\left(m_1^{\rho}\right)+\mathcal{G}\left(m_2^{\rho}\right)}{2}-\frac{\lambda \rho^{\lambda} \Gamma(\lambda+1)}{2\left(m_2^{\rho}-m_1^{\rho}\right)^{\lambda}}\left[^{\rho} I_{m_1^+}^{\lambda} \mathcal{G}\left(m_2^{\rho}\right)+^{\rho} I_{m_2^-}^{\lambda} \mathcal{G}\left(m_1^{\rho}\right)\right]\\&=\frac{m_2^{\rho}-m_1^{\rho}}{2} \int_{0}^{1}\left[\left(1-\zeta^{\rho}\right)^{\lambda}-\zeta^{\rho \lambda}\right] \zeta^{\rho-1} \mathcal{G}^{\prime}\left(\zeta^{\rho} m_1^{\rho}+\left(1-\zeta^{\rho}\right) m_2^{\rho}\right) d\zeta.
		\end{aligned}
		\end{equation}
	\end{footnotesize}
\end{lemma}

\begin{theorem}\label{ch2katthe}
	Suppose that  $\mathcal{G}: \left[m_1^{\rho}, m_2^{\rho}\right] \subset \mathbb{R}_{+} \rightarrow \mathbb{R}$ is a differentiable mapping on \(\left(m_1^{\rho}, m_2^{\rho}\right)\) with $0 \leq m_1<m_2$. If $|\mathcal{G}^{\prime}|$ is convex on \(\left[m_1^{\rho}, m_2^{\rho}\right],\) then the following inequality
	\begin{footnotesize}
		\begin{equation}
		\begin{aligned}
		&\left|\frac{\mathcal{G}\left(m_1^{\rho}\right)+\mathcal{G}\left(m_2^{\rho}\right)}{2}-\frac{\lambda \rho^{\lambda} \Gamma(\lambda+1)}{2\left(m_2^{\rho}-m_1^{\rho}\right)^{\lambda}}\left[^{\rho} I_{m_1^+}^{\lambda} \mathcal{G}\left(m_2^{\rho}\right)+^{\rho} I_{m_2^-}^{\lambda} \mathcal{G}\left(m_1^{\rho}\right)\right]\right|\\& \leq \frac{m_2^{\rho}-m_1^{\rho}}{2 \rho(\lambda+1)}\left(1-\frac{1}{2^{\lambda}}\right)\left[\left|\mathcal{G}^{\prime}\left(m_1^{\rho}\right)\right|+\left|\mathcal{G}^{\prime}\left(m_2^{\rho}\right)\right|\right]\label{ch2katineq}
		\end{aligned}
		\end{equation}
	\end{footnotesize}
	holds.
\end{theorem}

\begin{remark}
	Choosing $\rho=1$ in Theorem \ref{ch2katthe} reduces inequality (\ref{ch2katineq})
	to inequality (\ref{sfin}) in Theorem \ref{sfth}.
\end{remark}

Other important results involving Katugampola fractional integrals include the work of Mehreen and Anwar \cite{mehreen2018integral}, who generalized Theorem \ref{ths_Rch2} given as follows.

\begin{theorem}\label{mehrentheorem}
	Suppose that \(\lambda>0\) and \(\rho>0 .\) Let \(	\mathcal{G}:\left[m_1^{\rho}, m_2^{\rho}\right] \rightarrow \mathbb{R}\) be a positive function with \(0 \leq m_1<m_2\) and \(\mathcal{G} \in X_{c}^{p}\left(m_1^{\rho}, m_2^{\rho}\right).\) If \(	\mathcal{G}\) is also a convex function on \([m_1, m_2],\) then the following inequality
	\begin{footnotesize}
		\begin{equation}
		\begin{aligned}
		2^{s-1}\mathcal{G}\left(\frac{m_1^{\rho}+m_2^{\rho}}{2}\right)&\leq \frac{\rho^{\lambda} \Gamma(\lambda+1)}{2\left(m_2^{\rho}-m_1^{\rho}\right)^{\lambda}}\left[^{\rho} I_{m_1^+}^{\lambda} \mathcal{G}\left(m_2^{\rho}\right)+^{\rho} I_{m_2^-}^{\lambda} \mathcal{G}\left(m_1^{\rho}\right)\right]\\& \leq \bigg[\frac{\lambda}{(\lambda+s)}+\lambda\beta(\lambda,s+1)\bigg]\frac{\mathcal{G}\left(m_1^{\rho}\right)+\mathcal{G}\left(m_2^{\rho}\right)}{2}\label{mehreeninequality}\end{aligned}
		\end{equation}
	\end{footnotesize}
	holds.
\end{theorem}


Anderson \cite{anderson2016taylor} provided generalized H-H type inequalities involving conformable fractional integral as follows.
\begin{theorem}
Suppose that \(\alpha \in(0,1], m_{1}, m_{2} \in \mathbb{R}\) where \(m_{1}<m_{2}\), and \(\mathcal{G}:\left[m_{1}, m_{2}\right] \longrightarrow \mathbb{R}\) is an \(\alpha\) -
fractional differentiable function such that \(D_{\alpha}(\mathcal{G})\) is increasing. Then, we have

\begin{equation}
\frac{\alpha}{m_{2}^{\alpha}-m_{1}^{\alpha}} \int_{m_{1}}^{m_{2}} \mathcal{G}(y) \mathrm{d}_{\alpha} y \leq \frac{\mathcal{G}\left(m_{1}\right)+\mathcal{G}\left(m_{2}\right)}{2}\label{comforfir}
\end{equation}
In addition, if the mapping \(\mathcal{G}\) is decreasing on \(\left[m_{1}, m_{2}\right]\), then
\begin{equation}
\mathcal{G}\left(\frac{m_{1}+m_{2}}{2}\right) \leq \frac{\alpha}{m_{2}^{\alpha}-m_{1}^{\alpha}} \int_{m_{1}}^{m_{2}} \mathcal{G}(y) \mathrm{d}_{\alpha} y\label{comforse}
\end{equation}

\end{theorem}

\begin{remark}
	If \(\alpha=1\), then we clearly see that inequalities (\ref{comforfir}) and (\ref{comforse}) reduce to
	inequality (\ref{hh}).
\end{remark}


In \cite{set2018hermite} Set et.al. provided the generalized H-H type inequalities involving conformable fractional integral as follows.

\begin{theorem}\label{cth}
	Let \(\mathcal{G}:[m_1, m_2] \rightarrow \mathbb{R}\) be a function with \(0 \leq a<b\) and \(f \in L_{1}[m_1, m_2] .\) If \(\mathcal{G}\) is a convex on
	\([m_1, m_2],\), then the following inequality
	\begin{equation}
	\mathcal{G}\left(\frac{m_1+m_2}{2}\right) \leq \frac{\Gamma(\alpha+1)}{2(m_2-m_1)^{\alpha} \Gamma(\alpha-n)}\left[\left(I_{\alpha}^{m_1} \mathcal{G}\right)(m_2)+\left(^{m_2} I_{\alpha} \mathcal{G}\right)(m_1)\right] \leq \frac{\mathcal{G}(m_1)+\mathcal{G}(m_2)}{2}\label{cin}
	\end{equation}

holds, where \(\alpha \in(n, n+1]\).
\end{theorem}

Sarikaya et al. \cite{sarikaya2019hermite} present the following H-H inequalities for conformable fractional integral.

\begin{theorem}
	Let \(0<m_1<m_2, \alpha \in(0,1), \mathcal{G}:[m_1, m_2] \rightarrow \mathbb{R}\) be a convex function
	and \(J_{\alpha} f\) exists on \([m_1, m_2] .\) Then one has
	\begin{equation}
	\mathcal{G}\left(\frac{m_1^{\alpha}+m_2^{\alpha}}{2}\right) \leq \frac{\alpha}{m_2^{\alpha}-m_1^{\alpha}} \int_{m_1}^{m_2} \mathcal{G}\left(\gamma^{\alpha}\right) d_{\alpha} \gamma \leq \frac{\mathcal{G}\left(m_1^{\alpha}\right)+\mathcal{G}\left(m_2^{\alpha}\right)}{2}
	\end{equation}
\end{theorem}


For more results on generalization of H-H type inequalities, we refer the interested reader to  \cite{alomari2011some}, \cite{set2018certain}, \cite{set2018hermite} and \cite{agarwal2017some}.

\section{Applications to special means}\label{ch2application}
The following means for positive real numbers $m_1, m_2 \in \mathbb{R}_{+}$ $(m_1\neq m_2)$ exist in the literature \cite{bullen2003handbook,dragomir2004selected}.

\begin{itemize}
	\item [1.]The arithmetic mean:
	\begin{equation*}
	A(m_1, m_2)=\frac{m_1+m_2}{2}.
	\end{equation*}
	\item [2.]	The geometric mean:
	\begin{equation*}
	G(m_1, m_2)=\sqrt{m_1 m_2}.
	\end{equation*}
	\item [3.] The logarithmic mean:
	\begin{equation*}
	L(m_1, m_2)=\frac{m_2-m_1}{\ln m_2-\ln m_1}, m_1, m_2 \neq 0.
	\end{equation*}

	\item [4.]	Generalized log-mean:
	\begin{equation*}
	L_{\vartheta}(m_1, m_2)=\left[\frac{m_2^{\vartheta+1}-m_1^{\vartheta+1}}{(\vartheta+1)(m_2-m_1)}\right]^{1 / \vartheta}, \vartheta\in \mathbb{Z}\setminus\{-1,0\}.
	\end{equation*}
\end{itemize}

We note that \(L_{\vartheta}\) is monotonically increasing over \(\vartheta \in \mathbb{R}\) with \(L_{-1}=L\).
We, in particular, obtain the following inequality \(G \leq L \leq A\).
These special means can be frequently applied to numerical approximations, as well as other related problems that can be obtained in different fields. Several results that deal with the special means have been reported in the literature (see \cite{dragomir1998two} and \cite{sarikaya2015some}).

Dragomir and Agarwal \cite{dragomir1998applications} applying the results of Theorem\ref{ch2theorem1} to  establish the following new inequalities connecting the above means.

\begin{proposition}
 Let \(m_1, m_2 \in \mathbb{R}, m_1<m_2\) and \(\theta \in \mathbb{N}, \theta \geq 2 .\) Then, the following inequality
 \begin{equation*}
 \left|A\left(m_1^{\theta}, m_2^{\theta}\right)-L_{\theta}(m_1, m_2)\right| \leq \frac{\theta(m_2-m_1)}{4} A\left(|m_1|^{\theta-1},|m_2|^{\theta-1}\right)
 \end{equation*}
 holds.
\end{proposition}

\begin{proposition}
	Let \(m_1, m_2 \in \mathbb{R}, m_1<m_2\), and \(0 \notin[m_1, m_2] .\) Then, the following inequality
	\begin{equation*}
	\left|A\left(m_1^{-1}, m_2^{-1}\right)-\bar{L}^{-1}(m_1, m_2)\right| \leq \frac{(m_2-m_1)}{4} A\left(|m_1|^{-2},|m_2|^{-2}\right)
	\end{equation*}
	holds.
\end{proposition}

Furthermore, Kirmaci \cite{kirmaci2004inequalities}  established application to special means using the result of Theorem\ref{thmidclassch2} as follows.

\begin{proposition}
Let \(m_1, m_2 \in \mathbb{R}, m_1<m_2\) and \(\theta \in N, \theta \geqslant 2 .\) Then, we get
\begin{equation*}
\left|L_{\theta}(m_1, m_2)-A^{\theta}(m_1, m_2)\right| \leqslant \frac{\theta(m_2-m_1)}{4} A\left(|m_1|^{\theta-1},|m_2|^{\theta-1}\right)
\end{equation*}
\end{proposition}

 One can consult the following references \cite{bullen2013means,zhou2020new} for a comprehensive study on the special means.

\section{Applications to quadrature formula}

Let \(\mathcal{G}:[m_1, m_2] \rightarrow \mathbb{R}\) be a twice differentiable function on  $(m_1,m_2)$, such that \(\mathcal{G}^{\prime \prime}(x)\) is bounded on the given interval. This can be written as
\begin{equation*}
\left\|\mathcal{G}^{\prime \prime}\right\|_{\infty}=\sup _{x \in(m_1, m_2)}\left|\mathcal{G}^{\prime \prime}(x)\right|<\infty.
\end{equation*}
The following results are referred to as the midpoint and trapezoid inequalities,
\begin{equation}
\left|\int_{m_1}^{m_2} \mathcal{G}(x) d x-(m_2-m_1) \mathcal{G}\left(\frac{m_1+m_2}{2}\right)\right| \leq \frac{(m_2-m_1)^{3}}{24}\left\|\mathcal{G}^{\prime \prime}\right\|_{\infty}\label{ch2equa}
\end{equation}
and
\begin{equation}
\left|\int_{m_1}^{m_2} \mathcal{G}(x) d x-(m_2-m_1) \frac{\mathcal{G}(m_1)+\mathcal{G}(m_2)}{2}\right| \leq \frac{(m_2-m_1)^{3}}{12}\left\|\mathcal{G}^{\prime \prime}\right\|_{\infty},\label{ch2eq}
\end{equation}
respectively.

Therefore, the integral \(\int_{m_1}^{m_2} \mathcal{G}(x) d x\) can be approximated in terms of the midpoint formula
\begin{equation*}
\int_{m_1}^{m_2} \mathcal{G}(x) d x \cong(m_2-m_1) \mathcal{G}\left(\frac{m_1+m_1}{2}\right),
\end{equation*}
and the trapezoid formula
\begin{equation*}
\int_{m_1}^{m_2} \mathcal{G}(x) d x \cong(m_2-m_1) \frac{\mathcal{G}(m_1)+\mathcal{G}(m_2)}{2},
\end{equation*}
respectively.

The midpoint and trapezoid inequalities can be grouped in the most important relationship, the H-H inequality (\ref{hh}).

Suppose that \(d\) is a partition of the interval \([m_1, m_2]\) such that \(m_1=z_{0}<z_{1}<\cdots<z_{n-1}<z_{n}=m_2\). We write the following quadrature formula
\begin{equation*}
\int_{m_1}^{m_2} \mathcal{G}(x) d x=T_{i}(\mathcal{G}, d)+E_{i}(\mathcal{G}, d), i=1,2,
\end{equation*}
whereby
\begin{equation*}
T_{1}(\mathcal{G}, d)=\displaystyle{\sum_{i=0}^{n-1}} \frac{\mathcal{G}\left(z_{i}\right)+\mathcal{G}\left(z_{i+1}\right)}{2}\left(z_{i+1}-z_{i}\right)
\end{equation*}
is the trapezoidal version, and
\begin{equation*}
T_{2}(\mathcal{G}, d)=\displaystyle{\sum_{i=0}^{n-1}} \mathcal{G}\left(\frac{z_{i}+z_{i+1}}{2}\right)\left(z_{i+1}-z_{i}\right)
\end{equation*}
stands for the midpoint version.\\
The remainder term \(E_{1}(\mathcal{G}, d)\) for the integral $\int_{m_1}^{m_2} \mathcal{G}(x) d x$ estimated by the trapezoidal formula \(T_1(\mathcal{G}, d)\) satisfies
\begin{equation}
|E_1(\mathcal{G}, d)| \leq \frac{M}{12} \displaystyle{\sum_{i=0}^{n-1}}\left(z_{i+1}-z_{i}\right)^{3}.\label{Errorone}
\end{equation}
Meanwhile, that of the midpoint formula \(T_2(\mathcal{G}, d)\) satisfies
\begin{equation}
|E_2(\mathcal{G}, d)| \leq \frac{M}{24} \displaystyle{\sum_{i=0}^{n-1}}\left(z_{i+1}-z_{i}\right)^{3}.\label{Errortwo}
\end{equation}
These remainder terms (\ref{Errorone}) and (\ref{Errortwo})  can be used to estimate the error bounds of many numerical integration. Furthermore, the inequalities (\ref{ch2equa}) and (\ref{ch2eq}) can only hold if the second derivative is bounded on the interval $(m_1,m_2)$, and $\mathcal{G}$ is a twice differentiable function. This encourages many researchers to determine inequalities with less than or equal to one derivative. 

For example, Dragomir and Agarwal \cite{dragomir1998applications} estimated the reminder term through one derivative as follows.
\begin{proposition}
Let \(\mathcal{G}\) be a differentiable function on \(K^{o}, m_1, m_2 \in K^{o}\) with \(m_1<m_2.\) If \(\left|\mathcal{G}^{\prime}\right|\) is convex
on \([m_1, m_2]\), then the following
\begin{equation*}
\begin{aligned}|E(\mathcal{G}, d)| & \leq \frac{1}{8} \sum_{j=0}^{\theta-1}\left(z_{j+1}-z_{j}\right)^{2}\left(\left|\mathcal{G}^{\prime}\left(z_{j}\right)\right|+\left|\mathcal{G}^{\prime}\left(z_{j+1}\right)\right|\right) \\ & \leq \frac{\max \left\{\left|\mathcal{G}^{\prime}(m_1)\right|,\left|\mathcal{G}^{\prime}(m_2)\right|\right\}}{4} \sum_{j=0}^{n-1}\left(z_{j+1}-z_{j}\right)^{2} \end{aligned}
\end{equation*}
 holds.
\end{proposition}

Other important result was established by Kirmaci
\cite{kirmaci2004inequalities}, who estimated the reminder term through one derivative as follows.
\begin{proposition}
Let \(\mathcal{G}\) be a differentiable function on \(K^{o}, m_1, m_2 \in K^{o}\) with \(m_1<m_2.\) If \(\left|\mathcal{G}^{\prime}\right|\) is convex
on \([m_1, m_2]\), then the following
\begin{equation*}
|E(\mathcal{G}, d)| \leqslant \frac{1}{8} \sum_{j=0}^{\theta-1}\left(x_{j+1}-x_{j}\right)^{2}\left(\left|\mathcal{G}^{\prime}\left(x_{j}\right)\right|+\left|\mathcal{G}^{\prime}\left(x_{j+1}\right)\right|\right)
\end{equation*}
\end{proposition}

 Thus, these estimates remain as open-ended problems when considering their wider areas of applications \cite{dragomir1998two,kirmaci2004inequalities}.

\bibliographystyle{plain}
\bibliography{Review}
\end{document}